\documentclass[12pt,letterpaper]{article}
\usepackage{fullpage}
\usepackage[top=2cm, bottom=4.5cm, left=2.5cm, right=2.5cm]{geometry}
\usepackage{amsmath,amsthm,amsfonts,amssymb,amscd}
\usepackage{lastpage}
\usepackage{ulem}
\usepackage{enumerate}
\usepackage{pgfplots}
\pgfplotsset{compat=1.9}

\usepackage{algorithm} 
\usepackage{pgfplots}
\pgfplotsset{compat=1.5.1}
\usepackage{caption}
\usepackage{hyperref}
\usepackage{algpseudocode} 
\usepackage{fancyhdr}
\usepackage{mathrsfs}
\usepackage{listings}
\usepackage{graphicx}
\graphicspath{ {./images/} }

\newcommand\course{\text{\today}}
\newcommand\name{Michael Nisenzon} 

\begin{document}

\section*{Overview}
This is a summary of the splitting circle method as implemented$^1$ in Pari. The overall idea is to recursively split the polynomial $p(x)$ into two equal-degree factors $F(x), G(x)$ at each step until we reach linear factors which are the roots.

We can divide the algorithm into two parts: determining a splitting circle that contains roughly half of the roots on each side and estimating the factor $F(x)$ consisting of the product of the roots inside the circle. In this part, we will go in detail on using the splitting circle.

\section*{Finding a splitting circle}
\subsection*{Overview}
We start with a monic polynomial $P(x)= x^n + a_1x^{n-1} + ... + a_n$ of degree $n$ where the roots have moduli $\varrho_1(P),..., \varrho_n(P)$ arranged in increasing order. We would like to find a circle that contains some roots but not all and that the two collections of roots are as far apart as possible. Intuitively, if we see that the ratio $\frac{\varrho_{j+1}(P)}{\varrho_{j}(P)}$ is sufficiently large, then the geometric mean $\sqrt{\varrho_{j+1}(P)\varrho_{j}(P)}$ would make for a good radius for the splitting circle centered at the origin. Unfortunately, the roots of $P$ can have the exact same moduli.

Our algorithm has two steps: finding an appropriate center for the splitting circle and then calculating the radius.

\subsection*{Finding the center}
For our splitting circle, we care about the range of the moduli, as expressed by the ratio $\frac{\varrho_{n}(P)}{\varrho_{1}(P)}$. We would like to pick a center such that the range from the new center is sufficiently large.

First, we standardize the polynomial by shifting the coordinates to have the center of mass of the roots be the origin. Using Vieta's formulas, this is equivalent to considering the new polynomial $P_1(x) = P(x- a_1/n)$. Next, we would like to scale the coordinates such that the largest modulus is close to 1, which is the change of coordinates $P_2(x) = P_1(\varrho_n(P_1)x)$. Finally, for the shift, we consider four points outside the circle: $v = (2,2i,-2i,-2)$ and consider the shifted polynomials $Q_j = P_2(x-v_j)$. We choose as center the point $v_j$ that maximizes the ratio $\frac{\varrho_{n}(Q_j)}{\varrho_{1}(Q_j)}$. Let $Q$ be that maximizing polynomial. As the center of mass was defined to be the origin in the first step, we can easily show that $\frac{\varrho_{n}(Q)}{\varrho_{1}(Q)}$ is bounded below by a constant, specifically $e^{0.3}$. 

In the example below, we show the final shift of the center for a sample centered polynomial $P$. Note that the ratio $\frac{\varrho_{n}(Q)}{\varrho_{1}(Q)}$ substantially increases after choosing $(-2,0)$ as the new center.

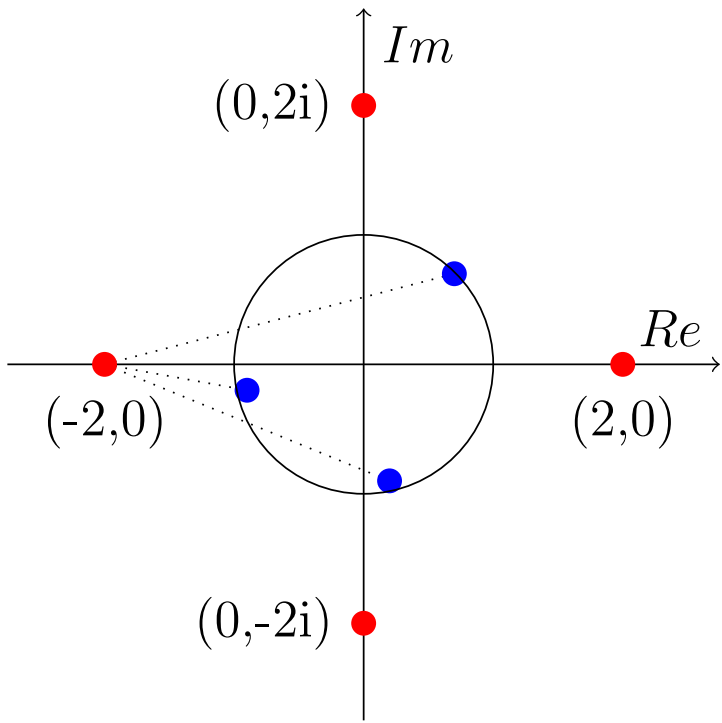
\begin{figure}[ht]
\resizebox{.6\textwidth}{!}{

\begin{tikzpicture}
\begin{axis}[ 
    ticks=none,
    axis lines = middle,
    axis line style={->},
    ymin=-5.5,ymax=5.5,
    xmin=-5.5, xmax=5.5,
    xlabel={$Re$},
    ylabel={$Im$},
    axis equal image
]
\addplot[dotted, domain=-4:-1.8] {-.181818*(x+4)};
\addplot[dotted, domain=-4:1.4] {0.259259259*(x+4)};
\addplot[dotted, domain=-4:0.4] {-0.4090909*(x+4)};

\node[label={180:{(0,2i)}},circle, red,fill,inner sep=2pt] at (axis cs:0,4) {};
\node[label={180:{(0,-2i)}},circle, red,fill,inner sep=2pt] at (axis cs:0,-4) {};
\node[label={270:{(2,0)}},circle, red,fill,inner sep=2pt] at (axis cs:4,0) {};
\node[circle, blue,fill,inner sep=2pt] at (axis cs:1.4,1.4) {};
\node[circle, blue,fill,inner sep=2pt] at (axis cs:-1.8,-.4) {};
\node[circle, blue,fill,inner sep=2pt] at (axis cs:.4,-1.8) {};
\node[label={270:{(-2,0)}},circle, red,fill,inner sep=2pt] at (axis cs:-4,0) {};
\draw (axis cs:0,0) circle [blue, radius=2];
\end{axis}
\end{tikzpicture}
}
\caption{Finalizing the choice of center of the splitting circle}

\end{figure}

\subsection*{Finding the radius}
After recentering, we now have a polynomial $P$ with ratio $\frac{\varrho_{n}(P)}{\varrho_{1}(P)} > e^{\Delta}$ for $\Delta > 0$. We would like to find an annulus $\Gamma$ that contains no roots of $P$ and such that $n/2$ roots are inside the interior circle.

We define $\Gamma$ by the interior radius $r$ and exterior radius $R$. We can then initialize $r_0 = \varrho_1(P)$ and $R= \varrho_n(P)$ and note that there is $i_0 = 1$ root in the inner circle and $j_0 = n-1$ roots in the outer circle. Afterwards, we will iteratively shrink our annulus until $n/2$ roots are located in the inner circle, with none in the annulus. At each step $\ell$, we take the geometric mean $\rho = \sqrt{r_{\ell-1}R_{\ell-1}}$ and choose our new annulus to be either $\Gamma_\ell = (r_{\ell-1}, \rho)$ or $\Gamma_\ell = (\rho, R_{\ell-1})$ using a subroutine that finds on the number of roots $k$ within the disk $\{|z| = \rho\}$. If $k < (i_{\ell-1} + j_{\ell-1})/2)$, we choose the former, otherwise we choose the latter. We stop when $i_\ell = j_\ell$, that is, when the annulus contains no roots. Then there must be an index $k$ such that $\varrho_k(P) < r_ell < R_\ell < \varrho_{k+1}(P)$. We conclude by finding $m= \varrho_k(P), M=\varrho_{k+1}(P)$ and setting our splitting circle to the geometric mean $\{|z| =  \sqrt{mM}\}$. 

We claim that the initial guarantee of $\frac{\varrho_{n}(P)}{\varrho_{1}(P)} > e^{\Delta}$ implies that $\frac{R_\ell}{r_\ell} > e^{\delta/(n-1)}$. Then scaling our splitting circle to a unit circle by a change of coordinates $z \rightarrow z/R$ gives us a root-free annulus of the form $(e^{-\delta}, e^{\delta})$ where $\delta = \frac{1}{2}\log(M/m)$.

\subsection*{Graeffe's method}
To numerically evaluate the modules of the roots of a polynomial $P$, we will use Graeffe's method, which follows from the observation that if 
\[P(x) = (x-z_1)(x-z_2)...(x-z_n)\] then 
\[P(x)P(-x) = (-1)^n(x^2-z_1^2)(x^2-z_2^2)...(x^2-z_n^2)\]
We note that the roots of $Q(x^2) := P(x)P(-x)$ are the squares of the roots of $P(x)$, and we will denote this operation as $Q = Graeffe(P)$. If 
\[P(x) = \sum_{i=0}^n a_i x^i\]

We can efficiently calculate $Graeffe(P)$ thanks to the formula
\[Q = A^2 - zB^2, \qquad A = \sum_i a_{2i} z^i, \qquad B = \sum_i a_{2i+1} z^i\]
We iteratively define $P_{m+1} = Graeffe(P_m)$ and note that 
\[P_m(x) = a_0^{(m)} + a_1^{(m)}x + ... + a_n^{(m)}x^n \]

Through Vieta's formulas, we see that at the $m$th iteration,
\[\frac{a_k^{(m)}}{a_n^{(m)}} \sum_{i_1 < ... < i_{n-k}} (z_{i_1} * ... * z_{i_{n-k}})^{2^m}\]
where $z_1,...,z_n$ are the roots of $P$. If the moduli are ordered $|z_1| < ... < |z_n|$, then we get the approximation
\[\frac{a_k^{(m)}}{a_0^{(m)}}=  (z_{k+1}*...*z_{n})^{2^m}\] so that
\[\varrho_k(P) = \lim_{m \rightarrow \infty} \left|\frac{a_k^{(m)}}{a_0^{(m)}}\right|^{2^{-m}}\]

Of course, using the Graeffe iteration directly only guarantees convergence when all roots have distinct moduli and there is no bound as to the rate of convergence. We will use Graeffe for three distinct cases: calculating the number of roots in a disk, calculating the kth modulus $\varrho_k(P)$, and calculating the $n$th modulus $\varrho_n(P)$, where there are more efficient algorithms for estimating the largest modulus than in the general case. 

\subsection*{Number of roots in a disk}
Given a polynomial $P$, error parameter $\tau$ and a radius $R$, we would like to find an index $k$ such that 
$\varrho_k(P)e^{-t} < R < \varrho_{k+1}(P)e^{t}$. Doing so tells us that there are exactly $k$ roots in the disk of radius $R$. Applying a change of coordinates, we can reduce this problem to the case where $R=1$. We will get an initial error tolerance of $t$ using a result of Schonhage and apply Graeffe's method until the $m$ th iteration of the error tolerance $t_m < \tau$. 

We will use the following theorem$^2$ from Schonhage without proof:

\textbf{Theorem 1} Let $P(x) = a_0 + a_1 x + ... + a_n x^n$ be a complex-valued polynomial with $k$ an integer such that $1\leq k \leq n$. Given the existence of $c,q > 0$ such that
\[|a_{k-m}| \leq cq^m|a_k|, \qquad m = 1:k\] we have that
\[\varrho_k(P) \leq (c+1)q(n-k+1)\]

\textbf{Corollary 1}
In addition, if we have that 
\[|a_{k+m}| \leq cq^m|a_k| , \qquad m = 1:n-k\]
then
\[\varrho_{k+1}(P) \geq \frac{1}{(c+1)q(n-k+1)}\]

\textit{Proof.} Consider the reciprocal polynomial $\tilde{P}(x) = x^nP(\frac{1}{x})$. We see that $\varrho_{n-k}(\tilde{P}) = \varrho_k(P)$ and we apply the inequality from the previous theorem.

Moreover, we see that $c=1,q=1$ satisfies the condition of the theorem for $k = \arg\max_{0 \leq i \leq n} |a_i|$. Then applying the theorem and corollary tells us that $\varrho_{k}(P) \leq 2n$ and $\varrho_{k+1}(P) \geq \frac{1}{2n}$. Setting $t > \log(2n)$ gives us the desired bound. However, this error tolerance is unworkably high, with many potentially valid values of $k$. As stated earlier, we can tighten the error bound with Graeffe's method.

Let $P_m(x) = b_0 + b_1(x) + ... + b_n x^n$ be the $m$th Graeffe iteration of $P$ and $k = \arg\max_{0 \leq i \leq n} |b_i|$. As before,
$\varrho_{k}(P_m) \leq 2n$ and $\varrho_{k+1}(P_m) \leq \frac{1}{2n}$. However, taking the $2^m$ roots of the moduli tells us that 
$\varrho_{k}(P) \leq e^{t_m}$ and $\varrho_{k+1}(P) \geq e^{-t_m}$ where $t_m = \frac{\log(2n)}{2^m}$. Given our fixed error tolerance $\tau$, we choose the smallest $m$ such that $t_m < \tau$, guaranteeing a sufficient separation of $\varrho_{k}(P) < 1 <  \varrho_{k+1}(P)$ up to multiplication by $e^{\tau}$.

Implementing this algorithm in double precision requires us to truncate coefficients of our polynomials. The following result$^3$ from Schatzle says that sufficiently small perturbations of the coefficients lead to small changes in the estimated $k$th modulus.

\textbf{Theorem 2} Let $P$ and $\hat{P}$ are complex polynomials with degree $n > 0$ such that $|\hat{P} - P| < \varepsilon|\hat{P}|$ and let $\mu > 1$ be an upper bound of $\varrho_k(P)$. If $\varepsilon < \mu^{-n}2^{-4n}$ then
\[|\varrho_k(\hat{P}) - \varrho_k(P)| \leq \frac{2\mu(1+\mu)\varepsilon^{1/n}}{1-4(1+\mu)\varepsilon^{1/n}}\]

From there, we get the following corollary:

\textbf{Corollary 2} Let $P$ and $\hat{P}$ be complex polynomials with degree $n > 0$, let $\tau > 0$ and $k$ an integer such that $0 \leq k \leq n$. If
\[|\hat{P} - P| < \varepsilon|\hat{P}|, \qquad  \varepsilon < 2^{-4n}\tau^ne^{-3n\tau/2}\] then
\[\varrho_k(\hat{P}) \geq e^{-3\tau/4} \implies \varrho_k(P) \geq e^{-\tau}\] and
\[\varrho_{k+1}(\hat{P}) \leq e^{3\tau/4} \implies \varrho_k(P) \leq e^{\tau}\]
so that running the algorithm on the truncated approximation for a radius of $R$ will find the largest modulus $\varrho_k < 1$ up to multiplication by $e^{-\tau}$.

\textit{Proof.}
To apply the theorem, set $\varepsilon \leq \mu^{-n}2^{-4n}$ for $\mu = \tau^{-1}e^{-3\tau/2}$, with $\tau < 1$. Then
\[|\varrho_k(\hat{P}) - \varrho_k(P)| \leq \frac{2\mu(1+\mu)\varepsilon^{1/n}}{1-4(1+\mu)\varepsilon^{1/n}} \leq  \]
\[\frac{2\mu(1+\mu)\mu^{-1}2^{-4}}{1-4(1+\mu)\mu^{-1}2^{-4}} =  \]
\[\frac{1+\mu}{8-2(1+\mu^{-1})} = \]
\[\frac{1+\mu}{6-2\mu^{-1}} = \]
\[\frac{1+\tau^{-1}e^{3\tau/2}}{6-2\tau e^{-3\tau/2}}\]

On the other hand, assuming the corollary doesn't hold, we see that 
\[|\varrho_k(\hat{P}) - \varrho_k(P)| > e^{-3\tau/4} - e^{-\tau}\]

Taking the difference of the two expressions  
\[\frac{d}{d\tau} \frac{1+\tau^{-1}e^{3\tau/2}+2e^{\tau/2}-2e^{3\tau/4}}{6-2\tau e^{-3\tau/2}} e^{-3\tau/4} - e^{-\tau} \] tells us that $e^{-3\tau/4} - e^{-\tau}$ is below our bound for $\tau <1$, so by nonnegativity of norms, $\varrho_k(P) \geq e^{-\tau}$.

The corresponding implication 
\[\varrho_{k+1}(\hat{P}) \leq e^{3\tau/4} \implies \varrho_k(P) \leq e^{\tau}\]
holds by Corollary 1.

We use this to create the subroutine NRD.

\begin{algorithm}[H]
	\caption{NRD: Calculates number of roots in a disk} 
	 \hspace*{\algorithmicindent} \textbf{Input}: $(P, R,\tau)$ where $P$ is a polynomial of degree $n \geq 2$, $R>0$ is the radius of the disk, and $\tau > 0$ is an error parameter \\
 \hspace*{\algorithmicindent} \textbf{Output}: An integer $k$ with $0 \leq k \leq n$ such that $\varrho_k(P)e^{-\tau} < R < \varrho_{k+1}(P)e^{\tau}$
	\begin{algorithmic}[1]
	\State Set $\tau_0 = \tau$ and $\varepsilon_0 = 2^{-4n}\tau_0^ne^{-3n\tau_0/2}$ and calculate the coefficients of the scaled polynomial $P_0 = P(Rz)$ with a relative precision of $\varepsilon_0/2(n+1)$.
	\State Round $P_0$ to a polynomial $\hat{P}_0$ such that $|P_0 - \hat{P}_0| < \frac{\varepsilon_0}{2}|\hat{P}_0|$.
	\State Set counter $m=1$.
	\While{$\frac{3}{4}\tau_{m-1} < \log(2n)$}
	\State Increment counter $m= m+1$.
	\State Set $P_m = Graeffe(\hat{P}_{m-1})$.
	\State Set $\tau_m = \frac{3}{4}\tau_{m-1}$ and $\varepsilon_0m = 2^{-4n}\tau_m^ne^{-3n\tau_m/2}$
	\State Round $P_m$ to a polynomial $\hat{P}_m$ such that $|P_m - \hat{P_m}| < \varepsilon_m|\hat{P}_m$
	\EndWhile
	\State With $\hat{P_m}(x) = b_0 + b_1x + ... + b_nx^n$, return the index $k = \arg\max_{0 \leq i \leq n} |b_i|$

	\end{algorithmic} 
\end{algorithm}
Let's consider the correctness of the algorithm. After the last step, we have the inequalities $\varrho_k(\hat{P}_m) < \exp(3\tau_m/4)$ and $\varrho_{k+1}(\hat{P}_m) > exp(-3\tau_m/4)$. By Corollary 2, we see that $\varrho_k(P_m)< e^{\tau_m}$ and $\varrho_{k+1}(P_m)> e^{-\tau_m}$. As $\varrho_j(P_i) = \varrho_j(\hat{P}_{i-1})^2$, we induct backwards using the definition of $\tau_m$ to see that $\varrho_k(\hat{P}_j) < e^{3\tau_j/4}$ and $\varrho_{k+1}(\hat{P}_j) > e^{-3\tau_j/4}$ for $0 \leq j \leq m-1$. In particular, for our choice of $k$, $\varrho_k(\hat{P_0}) < e^{3\tau/4}$ and  $\varrho_{k+1}(\hat{P_0}) > e^{-3\tau/4}$  so by Corollary 2, we conclude that \[\varrho_k(P(Rx))e^{-\tau} < 1 <\varrho_{k+1}(P(Rx))e^{\tau}\]
as desired.

\subsection*{Finding the kth modulus}
Next, we will use Graeffe's method to estimate the $k$th modulus $\varrho_k(P)$ to within a multiplicative error of $e^{\tau}$ for some error tolerance $\tau > 0$. Let $P(x) = a_0 + a_1x + ... + a_n x^n$. Then we will scale $P$ by some factor $\rho > 0$ with $\tilde{P}(x) = P(\rho x) = \tilde{a}_0 + \tilde{a}_1x + ... + \tilde{a}_n x^n$
such that there exist two integers $\ell, h$ with
\[\ell < k \leq h, \qquad |\tilde{a}_\ell| = |\tilde{a}_h|, |\tilde{a}_j| \leq |\tilde{a}_\ell| \quad \text{for } j=0,1,...,n\]
Then setting $c=q=1$ allows $\ell,h$ to satisfy Theorem 1, with
\[\frac{1}{2n} < \varrho_\ell(\tilde{P}) \leq \varrho_k(\tilde{P}) \leq \varrho_h(\tilde{P}) < 2n\]
We then apply Graeffe as before to tighten the bounds.

We can choose $\ell$ by considering the largest $x$-coordinate of the corners less than $k$ in the lower convex envelope of the points $M_j = (j, \log(|a_j|))$. Similarly, we can choose $h$ as the smallest $x$-coordinate of the corners greater than $k$ in the same envelope. We would like $|\tilde{a}_\ell| = |\tilde{a}_h|$, which requires $|a_{\ell}|\rho^{\ell}= |a_{h}| \rho^{h}$ or \[\rho = (\frac{|a_\ell|}{|a_h|})^{1/(h-\ell)}\] 
In addition, by convexity, as shown in the diagram below, the line passing through $\log(a_\ell)$ and $\log(a_h)$ upper bounds $\log(a_j)$ for $0 \leq j \leq n$. We switch to $\log(|a_j|)$ because $\log(\rho^x) = x\log(\rho)$ is monotone in $x$.  As a result, the linear upper bound remains an upper bound when composed with a monotone function $x 
\rightarrow x\log(\rho)$ for $0 \leq x \leq n$. Therefore, all points $\tilde{M}_j = (j, \log(|\tilde{a}_j|))$ are below the line connecting $\tilde{M}_\ell,\tilde{M}_h$. By our argument above, Theorem 1 is satisfied.

\begin{figure}[ht]
\resizebox{.6\textwidth}{!}{

\begin{tikzpicture}
\begin{axis}[ 
    axis lines = middle,
    axis line style={->},
    xtick={0,3,4,6,11},
    xticklabels={$0$,$\ell$, $k$, $h$, $n$},
    domain=0:2*pi
    ymin=-1.5,ymax=6.5,
    xmin=0, xmax=12,
    xlabel={$j$},
    ylabel={$\log(|a_j|)$},
    axis equal image
]
\addplot[black, domain=0:3] {1.16666*x-.5};
\addplot[black, domain=3:6] {.6666*x+1};
\addplot[dotted, domain=0:11] {.6666*x+1};
\addplot[black, domain=6:8] {-.5*x+8};
\addplot[black, domain=8:10] {-1*x+12};
\addplot[black, domain=10:11] {-1*x+12};

\node[circle, red,fill,inner sep=2pt] at (axis cs:0,-.5) {};
\node[circle, red,fill,inner sep=2pt] at (axis cs:1,-1) {};
\node[circle, red,fill,inner sep=2pt] at (axis cs:2,1) {};
\node[circle, red,fill,inner sep=2pt] at (axis cs:3,3) {};
\node[circle, red,fill,inner sep=2pt] at (axis cs:4,2) {};
\node[circle, red,fill,inner sep=2pt] at (axis cs:5,3.5) {};
\node[circle, red,fill,inner sep=2pt] at (axis cs:6,5) {};
\node[circle, red,fill,inner sep=2pt] at (axis cs:7,1) {};
\node[circle, red,fill,inner sep=2pt] at (axis cs:8,4) {};
\node[circle, red,fill,inner sep=2pt] at (axis cs:9,0) {};
\node[circle, red,fill,inner sep=2pt] at (axis cs:10,2) {};
\node[circle, red,fill,inner sep=2pt] at (axis cs:11,1) {};

\end{axis}
\end{tikzpicture}
}
\caption{Convex envelope of log(coefficients) }

\end{figure}
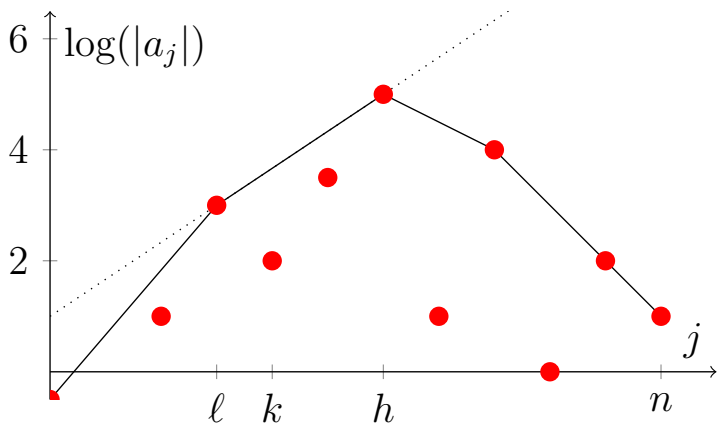

For convenience, we set $\rho$ to be a power of $2$, with \[\rho = 2^\beta, \qquad \beta = \left\lfloor \frac{1}{h-\ell}\log_2 \frac{|a_\ell|}{|a_h|}+1/2\right\rfloor\]

We then see that for $m \geq 0$, 
\[|\tilde{a}_{h-m}| \leq 2^{m/2}|\tilde{a}_h|, \qquad |\tilde{a}_{\ell+m}| \leq 2^{m/2}|\tilde{a}_\ell| \]
so the conditions of Theorem 1 are valid for $c=1,q=\sqrt{2}$. Due to rounding error, it is better to work with an upper bound of $c=1,q=3/2$, which gives us looser initial bounds of 
\[\frac{1}{3n} < \varrho_\ell(\tilde{P}) \leq \varrho_k(\tilde{P}) \leq \varrho_h(\tilde{P}) < 3n\]

We then use the following result$^1$ from Gourdon to guarantee the accuracy of our estimates given a truncated initial polynomial.

\textbf{Theorem 3} Let $P, \hat{P}$ be two polynomials of degree $n >0$, let $\tau >0$ and $k$ an integer such that $0 \leq k \leq n$. If
\[|\hat{P} - P| < \varepsilon|\hat{P}|, \qquad  \varepsilon < 2^{-n-1}(3n)^{-n}\tau^ne^{-3n\tau/2} \qquad \frac{1}{3n} < \varrho_k(P) < 3n\]
Then
\[\varrho_k(\hat{P})e^{-\tau} \leq \varrho_k(P) \leq \varrho_k(\hat{P})e^\tau\]
\begin{algorithm}[H]
	\caption{MOD: Evaluates the moduli of the roots.} 
	 \hspace*{\algorithmicindent} \textbf{Input}: $(R,k,\tau)$  where $P(x) = a_0 + a_1x + ... + a_nx^n$ is a polynomial of degree $n \geq 2$, $k$ is an integer with $1 \leq k \leq n$, and $\tau > 0$ is an error parameter.\\
 \hspace*{\algorithmicindent} \textbf{Output}: A floating point number $R>0$, such that $Re^{-\tau} \leq \varrho_k(P) \leq Re^\tau$.
	\begin{algorithmic}[1]
	\State If $a_0 = a_1 = ... = a_{n-1}=0$, return $R = 0$.
	\State Otherwise, set $P_0 = P$ and $\tau_0 = \tau/8$.
	\State Calculate $\tilde{P}_0(x) = P_0(\rho_0x)$ for $\rho_0 := 2^\beta, \quad \beta := \left\lfloor \frac{1}{h-\ell}\log_2 \frac{|a_\ell|}{|a_h|}+1/2\right\rfloor$ where $h,\ell$ are defined as in the convex envelope for $P_0$.
	\State Set $M$ to be the smallest natural number such that $2^{-M}\log(3n) < \tau/2$
	\For{$m=1,2,...,M$}
	\State Round the polynomial $\tilde{P}_{m-1}$ to a polynomial $\hat{P}_{m-1}$ such that \[|\tilde{P}_{m-1} - \hat{P}_{m-1}| < \varepsilon_{m-1}, \quad \varepsilon_{m-1} := 2^{-(n+1)}(3n)^{-n}\tau_{m-1}^ne^{-3n\tau_{m-1}/2}\]
	\State Calculate $P_m = Graeffe(\hat{P}_{m-1})$
	\State Calculate $\tilde{P}_m(x) = P_m(\rho_m x)$ where $\rho_m$ is defined as in the convex envelope for $P_m$.
	\State Set $\tau_m = \frac{3}{2}\tau_{m-1}$
	
	\EndFor
	\State Return $R = \rho_0\rho_1^{2^{-1}}\rho_2^{2^{-2}}...\rho_M^{2^{-M}}$
	\end{algorithmic} 
\end{algorithm}

Let us evaluate the correctness of our algorithm. For each step $m$, by our change of coordinates we get $\varrho_k(P_m) = \rho_m\varrho_k(\tilde{P}_m)$. After applying the perturbation result of Theorem 3, we see that
\[\varrho_m\varrho_k(\hat{P}_m)e^{-\tau} < \varrho_k(P_m) < \varrho_m\varrho_k(\hat{P}_m)e^{\tau}\]
Using the relation $\varrho_k(P_{m+1}) = \varrho_k(\hat{P}_m)^2$ gives us
\[\rho_0...\rho_M^{2^{-M}}\varrho_k(\tilde{P}_M)^{2^{-M}}e^{-\tau'} < \varrho_k(P_0) < \rho_0...\rho_M^{2^{-M}}\varrho_k(\tilde{P}_M)^{2^{-M}}e^{\tau'}, \quad \tau' = \tau_0 + \frac{\tau_1}{2} + ... + \frac{\tau_{M-1}}{2^{M-1}}\]

From Theorem 3, we see that as $\frac{1}{3n} < \varrho_k(\tilde{P}_M) < 3n$, then 
\[e^{-\tau/2} < \varrho_k(\tilde{P}_M)^{2^{-M}} < e^{\tau/2}\]
Moreover, \[0 \leq \tau' = \frac{\tau}{8}\left(1+ \frac{3}{4} + \left(\frac{3}{4}\right)^2 + ...\right) \leq \frac{\tau}{2}\]
so we conclude that
\[Re^{-\tau} < \varrho_k(P_0) < Re^{\tau}\]
for $R = \rho_0...\rho_M^{2^{-M}}$.
\subsection*{Upper bounding the moduli of the roots}
For the sketch of our splitting circle algorithm, we need to calculate $\varrho_n(P)$ often. It turns out that this can be accomplished more efficiently than $\varrho_k(P)$ for a general $k$. Again, given an error parameter $\tau>0$, we would like to find an $R > 0$ equal to $\varrho_n(P)$ up to multiplicative error,
such that 
\[Re^{-\tau} \leq \varrho_n(P) \leq Re^{\tau}\]
The idea of the algorithm is the same as in the previous case. However, we would like to use a different scalar $\rho$. Given an initial $P = a_0 + a_1x + ... + a_nx^n$, we define $\tilde{P}(x) = P(\rho x)$ such that
\[\frac{|\tilde{a}_j|}{|\tilde{a}_0|} \leq 2^j {n \choose j} \quad \text{for } j = 1,...,n \quad \text{and there exists } h, 1 \leq h \leq n \text{ such that } \frac{|\tilde{a}_j|}{|\tilde{a}_0|} \geq {n \choose h}\]
and $\rho$ can be determined by \[\rho = 2^\beta, \qquad \beta = \max_{1 \leq j \leq n} \left\lfloor \frac{1}{j}\log \frac{|a_\ell|}{|a_0|{n \choose j}}\right\rfloor\]

By Vieta's equations, we get the inequality
\[|\tilde{a}_h| \leq |\tilde{a}_0|{ n \choose h} \varrho_n(\tilde{P})^h\] so $\varrho_n(\tilde{P}) \geq 1$ when we account for the second assumption for our choice of $\rho$. Setting $c = 1, q=2n$, we see that $2^j{n \choose j} \leq (2n)^j$, so Theorem 1's conditions are satisfied and give an upper bound of $1 \leq \varrho_n(\tilde{P}) \leq 4n$. We continue with a perturbation result.

\textbf{Theorem 4} Let $P,\hat{P}$ be polynomials of degree $n > 0$ such that 
\[|P - \hat{P}| < |\beta|e^{-n\tau}\tau^n\]
where $\beta$ is the leading coefficient of $P$. Then if $\varrho_n(P) \geq 1$, we have that
\[\varrho_n(\hat{P})e^{-\tau} < \varrho_n(P) < \varrho_n(\hat{P})e^\tau\]

\textit{Proof.} First, let $\varrho_n(\hat{P}) > \varrho_n(P)$ and let $z$ designate a root of $\hat{P}$ with modulus $\varrho_n(P) > 1$. Then we have that
\[|P(z)| = |P(z) - \hat{P}(z)| \leq |z|^n|P-\hat{P}| \leq \varrho_n(\hat{P})^n|\beta|e^{-n\tau}\tau^n\]
Moreover, we have
\[|P(z)| \geq |\beta|\prod_{i=1}^n (|z|- \varrho_i(P)) \geq |\beta|(\varrho_n(\hat{P})-\varrho_n(P))^n\]

which gives us the bound 
\[\varrho_n(\hat{P})-\varrho_n(P) \leq \varrho_n(\hat{P})e^{-\tau}\tau\]

Applying the inequality $1+x\leq e^x$ gives us the desired bounds. An analogous argument holds for $\varrho(\hat{P}) < \varrho_n(P)$ by switching the roles of the two terms.

We now describe the algorithm to find the largest modulus of the roots. The code and justification follow largely from the previous example. 

\begin{algorithm}[H]
	\caption{MODMAX: Evaluates the largest modulus of the roots.} 
	 \hspace*{\algorithmicindent} \textbf{Input}: $(R,\tau)$  where $P(x) = a_0 + a_1x + ... + a_nx^n$ is a polynomial of degree $n \geq 2$ and $\tau > 0$ is an error parameter.\\
 \hspace*{\algorithmicindent} \textbf{Output}: A floating point number $R>0$, such that $Re^{-\tau} \leq \varrho_n(P) \leq Re^\tau$.
	\begin{algorithmic}[1]
	\State If $a_0 = a_1 = ... = a_{n-1}=0$, return $R = 0$.
	\State Otherwise, set $P_0 = P$ and $\tau_0 = \tau/8$.
	\State Calculate $\tilde{P}_0(x) = P_0(\rho_0x)$ for $\rho = 2^\beta, \qquad \beta = \max_{1 \leq j \leq n} \left\lfloor \frac{1}{j}\log \frac{|a_\ell|}{|a_0|{n \choose j}}\right\rfloor$
	\State Set $M$ to be the smallest natural number such that $2^{-M}\log(4n) < \tau/2$
	\For{$m=1,2,...,M$}
	\State Round the polynomial $\tilde{P}_{m-1}$ to a polynomial $\hat{P}_{m-1}$ such that \[|\tilde{P}_{m-1} - \hat{P}_{m-1}| < \varepsilon_{m-1}, \quad \varepsilon_{m-1} := |\beta+{m-1}\tau_{m-1}^ne^{-n\tau_{m-1}}\] \State where $\beta_{m-1}$ is the leading coefficient of $\tilde{P}_{m-1}$.
	\State Calculate $P_m = Graeffe(\hat{P}_{m-1})$
	\State Calculate $\tilde{P}_m(x) = P_m(\rho_m x)$ where $\rho_m$ is defined in terms of $P_m$ for $P_0$.
	\State Set $\tau_m = \frac{3}{2}\tau_{m-1}$
	
	\EndFor
	\State Return $R = \rho_0\rho_1^{2^{-1}}\rho_2^{2^{-2}}...\rho_M^{2^{-M}}$.
	\end{algorithmic} 
\end{algorithm}

Applying the algorithm to the reciprocal $P^*(x) = x^nP(\frac{1}{x})$ allows us to also find the smallest modulus of the roots.

\begin{algorithm}[H]
	\caption{MODMIN: Evaluates the smallest modulus of the roots.} 
	 \hspace*{\algorithmicindent} \textbf{Input}: $(R,\tau)$  where $P(x) = a_0 + a_1x + ... + a_nx^n$ is a polynomial of degree $n \geq 2$ and $\tau > 0$ is an error parameter.\\
 \hspace*{\algorithmicindent} \textbf{Output}: A floating point number $R>0$, such that $Re^{-\tau} \leq \varrho_1(P) \leq Re^\tau$.
	\begin{algorithmic}[1]
	\State If $a_0 = 0$, return $R = 0$.
	\State Otherwise, define the reciprocal polynomial $Q(x) = a_n + a_{n-1}(x) + ... + a_0x^n$ of $P$ and return 1/MODMAX($Q,\tau$).
	\end{algorithmic} 
\end{algorithm}

We now have all the tools we need to find the splitting circle of a polynomial.

\subsection*{Algorithms}
We can use the previous algorithms to to construct a new algorithm to find a splitting circle. We will break the search into two parts: finding the center and finding the radius of the splitting circle.

First, assume the center has already been determined. Finding the center should also give us an annulus $r < |x| < R$ along with two indices $i,j$ such that 
\[1\leq i < j \leq n-1, \qquad \varrho_i(P) < r < R < \varrho_{j+1}(P)\]
We now have to find a radius $\rho$ in the interval $(r,R)$ from which we can factorize $P$ into two nontrivial factors. The general idea was discussed above, and we will deal with practical implementation here.

\begin{algorithm}[H]
	\caption{RAD: Returns the radius of the splitting circle.} 
	 \hspace*{\algorithmicindent} \textbf{Input}: $(P, r, R, i, j)$  where $P$ is a polynomial of degree $n \geq 2$ that satisfies the centering conditions above\\
 \hspace*{\algorithmicindent} \textbf{Output}: $(\rho,k,\delta)$, where $\rho > 0$ is such that the circle $|z| = \rho$ contains $k$ roots of $P$ with $i \leq k \leq j$ and $\delta >0$ such that no root of $P$ is within the annulus $\rho e^{-\delta} < |z| < \rho e^{\delta}$
	\begin{algorithmic}[1]
	\State We calculate the error margin $\Delta := \log(R/r)$
	\item If $i = j$, then the annulus $r < |z| < R$ doesn't contain any roots of $P$ and we can do the following steps:
	\State \qquad We approximate the moduli $\varrho_i(P)$ and $\varrho_{i+1}(P)$ by 
	
	\quad $m= MOD(P,i,\Delta/8)$ and $M= MOD(P,i+i,\Delta/8)$
	\State \qquad We return $(\rho,\delta)$ where $\rho = \sqrt{mM}$ and $\delta = \frac{1}{2}\log(M/m) - \frac{\Delta}{8}$.
	\State Otherwise, we set $\rho = \sqrt{rR}, \delta = \frac{\Delta}{8(j-i)}$ and $k = NRD(P, \rho, \delta)$
	\State \qquad If $k < (i+j)/2$, return $RAD(R,r, \rho e^{-\delta},i,k)$.
	\State \qquad If $k > (i+j)/2$, return $RAD(R,\rho e^{\delta},R,k,j)$.
	\State \qquad If $k = (i+j)/2$ and $k<n/2$, return $RAD(R,r, \rho e^{-\delta},i,k)$. 
	\State Otherwise if $k = (i+j)/2$ and $k\geq n/2$, return $RAD(R,\rho e^{\delta},R,k,j)$.
	
	\end{algorithmic} 
\end{algorithm}

We see that $\delta$ is always greater than $c/(j-i+1)$ for some positive constant $c$. Furthermore, when the condition $i=j$ is reached, then we have an annulus $r<|z| < R$ that contains no roots of $P$ and so our final radius does give us a splitting circle.

In later sections, we will come up with an algorithm $FCS(P, k, \delta, \varepsilon)$ that splits $P$ into two approximate factors $\hat{F},\hat{G}$ such that $|P - \hat{F}\hat{G}| < \varepsilon |P|$ given a unit splitting circle. Combining $FCS$ with $RAD$ gives us an algorithm that factorizes a centered polynomial $P$.

\begin{algorithm}[H]
	\caption{HOM: Factorizes from a dilation of the splitting circle} 
	 \hspace*{\algorithmicindent} \textbf{Input}: $(P, r, R, i, j, \varepsilon)$  where $P$ is a polynomial of degree $n \geq 2$ that satisfies the centering conditions above, $\varepsilon > 0$ is an error parameter\\
 \hspace*{\algorithmicindent} \textbf{Output}: $(\hat{F},\hat{G})$, which are approximations of factors $F,G$ such that $|P-\hat{F}\hat{G}| < \varepsilon|P|$
	\begin{algorithmic}[1]
	\State Find $(\rho,k,\delta) = RAD(P,r,R,i,j)$.
	\State Find an approximation $\hat{Q}$ of $Q(z) = P(\rho z)$ with a relative precision of $\varepsilon'/n$ for $\varepsilon' = \frac{1}{4}\min(\rho^{-n},\rho^n)\varepsilon$
	\State Calculate $(F_0,G_0) = FCS(Q,k,\delta, \varepsilon')$
	\State Rescale the factors $\hat{F}(z) = F_0(z/\rho), \hat{G}(z) = G_0(z/\rho)$ with a relative precision of $\varepsilon'' = 2^{-(n+4)}\varepsilon/n$ and return the scaled factors
	$(\hat{F},\hat{G})$.
	\end{algorithmic} 
\end{algorithm}
We will now prove that the precision specified gives us $|P-\hat{F}\hat{G}| < \varepsilon|P|$. In the second step, the approximation $\hat{Q}$ of $Q$ satisfies 
\[|\hat{Q}-Q| < \varepsilon'|Q|\]
Moreover, $F_0,G_0$ satisfy
\[|\hat{Q}-\hat{F}\hat{G}| < \varepsilon'|\hat{Q}| < \varepsilon'|\hat{Q}| < 2\varepsilon'|Q|\] so 
\[|Q-\hat{F}\hat{G}| < 3\varepsilon'|Q|\]
The dilation scales errors by a factor at most $\max(\rho^n,\rho^{-n})$ so scaling tells us that
\[|P-FG| < \frac{3}{4}\varepsilon|P|\]
where $F(z) = F_0(z/\rho), G(z) = G_0(z/\rho)$.
Next, from our choice of $\varepsilon''$, we have the inequalities
$|F-\hat{F}|<2^{-n-4}\varepsilon|F|$ and $|G-\hat{G}|<2^{-n-4}\varepsilon|G|$. Multiplying these together tells us that 
\[|FG- \hat{F}\hat{G}| < 2^{-n-4}\varepsilon(|F|\cdot |G| + |F| \cdot |\hat{G}|) \leq 2^{-n-2}\varepsilon(|F|\cdot |G|)\]

From Schonhage, we get the inequality $|F|\cdot |G| \leq 2^{n-1}|FG|$, which is slightly tighter than the classical inequality $|F|\cdot |G| \leq 2^{n}|FG|$. This leads to the bound
\[|FG- \hat{F}\hat{G}| < \frac{\varepsilon}{8}|FG| \leq \frac{\varepsilon}{4}P\]
Putting all of our inequalities together,
\[|P - \hat{F}\hat{G}| \leq |P-FG| + |FG- \hat{F}\hat{G}| < \varepsilon|P|\]
as desired.

We now proceed to find the center of the splitting circle. The center $u$ for the polynomial $Q(z) = P(z+u)$ has to be chosen so that the modular ratio $\frac{\varrho_n(Q)}{\varrho_1(Q)}$ is sufficiently large. In practice, we need to preprocess to assure the calculations are well-conditioned. This is realized by our algorithm $CTR0$ which calls on the centering algorithm $CTR$ proper.

\begin{algorithm}[H]
	\caption{CTR0: Factorizes a given polynomial of degree $n \geq 2$} 
	 \hspace*{\algorithmicindent} \textbf{Input}: $(P, \varepsilon)$  where $P$ is a polynomial of degree $n \geq 2$, $\varepsilon > 0$ is an error parameter\\
 \hspace*{\algorithmicindent} \textbf{Output}: $(\hat{F},\hat{G})$, which are approximations of factors $F,G$ such that $|P-\hat{F}\hat{G}| < \varepsilon|P|$
	\begin{algorithmic}[1]
	\State If the constant coefficient $a_0$ of $P$ satisfies $|a_0| < \varepsilon|P|$, return $(z, (P(z)-P(0))/z)$.
	\State Otherwise, set $k_1 = NRD(P,1.9,0.05)$. If $k_1=n$, then $\varrho_n(P) \leq 1.9e^{.05}<2$ and we return $CTR(P,\varepsilon)$.
	\State Otherwise, consider the reciprocal polynomial $Q = P^* z^nP(1/z)$ and calculate $k_2 = NRD(Q,1.9,0.05)$. If $k_2=n$, then calculate $(F_0,G_0) = CTR(Q,\varepsilon)$ and return $(F_0^*,G_0^*)$.
	\State Otherwise, set $R = 1.9e^{0.05}$, $r = 1/R$, and return $HOM(P,r,R,n-k_2,k_1-1,\varepsilon)$.
	\end{algorithmic} 
\end{algorithm}

The algorithm $CTR0$ allows us to use $CTR$, which requires that the moduli of the roots are not too high or even directly move to $HOM$ if the constraints are satisfied, giving us a stability condition.

\begin{algorithm}[H]
	\caption{CTR: Finds the center of the splitting circle} 
	 \hspace*{\algorithmicindent} \textbf{Input}: $(P, \varepsilon)$  where $P(z)=a_0+a_1z+...+a_nz^n$ is a polynomial of degree $n \geq 2$ with $\varrho_n(P) \leq 2$, $\varepsilon > 0$ is an error parameter\\
 \hspace*{\algorithmicindent} \textbf{Output}: $(\hat{F},\hat{G})$, which are approximations of factors $F,G$ such that $|P-\hat{F}\hat{G}| < \varepsilon|P|$
	\begin{algorithmic}[1]
	\State Set $P_0(z) = P(z+u)$ where $u = \frac{-a_1}{n a_0}$ is the center of mass of the roots.
	\State Find an approximation $\hat{P}_0$ of $P_0$ with a relative precision of $\varepsilon_0' = \frac{1}{4}\varepsilon(1+|u|)^{-n}$. 
	\State Set $\varepsilon_0 = \frac{\varepsilon}{4}(1+|u|)^{-n}\frac{|P|}{|P_0|}$
	\State If the constant term $c$ of $\hat{P}_0$ satisfies $|c| < \varepsilon_0|\hat{P}_0|$ then skip to the last step with the polynomials $\hat{F}_0(z) = z, \hat{G}_0(z) = (\hat{P}_0(z)-c)/z$.
	\State Calculate an approximation $r = MODMAX(\hat{P}_0,0.01)$ of the largest modulus of the roots of $\hat{P}_0$, then check that $\rho = re^{-0.01} < 4$.
	\State Calculate an approximation $\hat{P}_1$ of the dilation $P_1(z) = \rho^{-n}\hat{P}_0(\rho z)$ with a relative precision of $\varepsilon_1/n$ for $\varepsilon_1 = \frac{1}{4}\varepsilon_0\sup(1,\rho)^{-n}$
	\State For $j=0,1,2,3$ calculate the polynomial $Q_j = \hat{P}_1(z+2e^{ij\pi/2})$ and the values $R_j = MODMAX(Q_j,0.01), r_j = MODMIN(Q_j,0.01)$.
	\State Let $j_0$ be the index that maximizes the ratio $R_j/r_j$. Set $P_2 = Q_{j_0}$ and $v = 2e^{ij_0\pi/2}$.
	\State Calculate $(F_2,G_2) = HOM(P_2, r,R, 1, n-1, \varepsilon_2)$ with $r= r_{j_n}e^{0.01}$, $R=R_{j_n}e^{-0.01}$ and $\varepsilon_2 = \varepsilon_1 \frac{|\hat{P}_0|}{|P_2|}3^{-n}$.
	\State Calculate $F_1(z) = F_2(z-v), G_1(z) = G_2(z-v)$.
	\State Find the approximations $\hat{F}_0, \hat{G}_0$ of the scaled polynomials $F_0 = \rho^kF_1(z/p), F_0 = \rho^{n-k}G_1(z/p)$ where $k=deg(F_1)$. Work with a relatie precision of $\varepsilon'/n$ for $\varepsilon' =2^{-n-4}\varepsilon_0$.
	\State Calculate the approximations $\hat{F},\hat{G}$ of $F = \hat{F}_0(z-u), G = \hat{G}_0(z-u)$ to a relative accuracy of $\varepsilon''/n$ for $\varepsilon = 2^{-n-4}\varepsilon$ and return $(\hat{F},\hat{G})$.
	\end{algorithmic} 
\end{algorithm}
We will now show that the precision of our calculations is sufficient. After calling HOM in line 9, we have that $|P_2 - F_2G_2| < \varepsilon_2|P_2|$, and the transformation $R(z) \rightarrow R(z-v)$ could amplify errors by a factor of at most $(1+|v|)^n = 3^n$ at most, giving us a bound of
\[|\hat{P}_1 - F_1G_1| < 3^n\varepsilon_2|P_2| = \varepsilon_1|\hat{P}_0|\]
As $|\hat{P}_1 - P_1| < \varepsilon_1 |\hat{P}_0|$, we see that
\[|P_1 - F_1G_1| < 2\varepsilon_1|\hat{P}_0|\]
Meanwhile, the dilation in line 11 gives us
\[|\hat{P}_0 - F_0G_0| < 2\varepsilon_1\sup_{1,\rho}^n|\hat{P}_0| = \frac{\varepsilon_0}{2}|\hat{P}_0|\]

We repeat the argument to show that 
the precision used at step 11 to find $\hat{F}_0,\hat{G}_0$ implies
\[|\hat{P}_0 - \hat{F}_0\hat{G}_0| < \varepsilon_0|\hat{P}_0|\]
Next, as $|\hat{P}_0 - P_0| < \varepsilon_0|\hat{P}_0|$ we see by the triangle inequality that
\[|P_0 - \hat{F}_0\hat{G}_0| < 2\varepsilon_0|\hat{P}_0|\] and therefore
\[|P_0 - FG| < 2\varepsilon_0(1+|u|)^n|\hat{P}_0|= \frac{\varepsilon}{2}|P|\]
The precision used to calculate $\hat{F},\hat{G}$ ensures that $|P- \hat{F}\hat{G}| < \varepsilon|P|$, which guarantees the desired output.

Note that for $|u| \leq 2, \rho \leq 4$, we have that 
\[\varepsilon_2 =\frac{1}{4}(1+|u|)^{-n}\sup(1,\rho)^{-n}3^{-n}\frac{|P|}{|P_2|}\varepsilon \geq \frac{1}{4}36{-n}\frac{|P|}{|P_2|}\varepsilon\]
In addition, we have that
\[|\hat{P}_0| \leq 2|P_0| \leq 2(1+|u|)^n|P| \leq 2\cdot 3^n|P|\]
Next, all the roots of $P_1$ are in the interior of the unit disk, so if $\alpha$ is the leading coefficient of $P_1$, then $|P_1| \leq |\alpha|2^n \leq |\hat{P}_0|2^n$ as $\alpha$ is also the dominant coefficient of $\hat{P}_0$. We then have that  $|\hat{P}_1| \leq |\hat{P}_0|2^{n+1}$. Using the inequality $|P_2| \leq 3^n|\hat{P}_1|$ we get that
\[\frac{|P_2|}{|P|}\leq 3^{2n}\cdot 2^{n+2} = 4*18^n\]
Putting all of our inequalities alows us to conclude that $\varepsilon_2 \geq \frac{1}{16}648^{-n}\varepsilon$. This lower bound of $\varepsilon_2$ makes sure that it is sufficiently large with respect to $\varepsilon$ that the calculations are well-conditioned.

\section*{Using a splitting circle to find the factors}
\subsection*{Overview}
Assume that we have a splitting circle $C$ that contains roughly half of the roots of our polynomial $p(x)$. Then by a change of coordinates, we can assume that our splitting circle is the unit circle. We would like to find the factor $F(x)$ that consists of the product of the roots in the circle. 

We first find an initial approximation $F_0$ of $F$ using the residue theorem. By polynomial division, we get the approximation $G_0$ of the remainder. We then use a form of Newton's method to get better approximations $F_1,G_1$ of $F,G$. We repeat the process until $F_1,G_1$ are linear factors.

\subsection*{Finding $F_0$}
By the residue theorem, we note that along the boundary of our unit circle $C$:
\[\frac{1}{2i\pi}\oint_C \frac{p'(z)}{p(z)}z^m dz = \sum_{|z| < 1, p(z) = 0}z^m = W_m\]
Using Newton's identities, we can turn the sums of powers from $m=1,...,n$ to the coefficients $\phi_i$ of $F$ using the recursive formulas:
\[\phi_1 = -W_1\]
\[\phi_m = -\frac{1}{m}(W_1\phi_{m-1} + ... + W_{m-1}\phi_1 + W_m)\]

For a large enough $N$, we can approximate the contour integral $W_m$ with the discrete sum
\[W_m = \frac{1}{N}\sum_{j=0}^{N-1} \frac{p'(\omega^j)}{p(\omega^j)}\omega^{(m+1)j}\]
where $\omega$ is the $N$th root of unity. Applying Newton's identities gives us our estimates $\phi_{i}$ and our initial approximation $F_0$ of the factor $F$, which should be within $O(e^{-\delta N})$ of the true value.

\subsection*{Newton-Schonhage iteration}
For our final approximation, we would like to find correction factors $f,g$ with degrees $deg(f) < deg(F_0), deg(g) < deg(G_0)$ such that $F_1 := F_0 + f, G_1 := G_0 + g$ are sufficiently good approximations of the true factors $F,G$. Specifically, we would like to minimize the L1 norm
\[\min_{f,g} ||P - (F_0 + f)(G_0 + g)||_1\] To the first order, we have that
\[P - (F_0 + f)(G_0 + g) \approx P - F_0G_0 - fG_0 - gF_0\]
so we would like to find $f,g$ such that 
\[P - F_0G_0 = fG_0 + gF_0\]
We can ignore $g$ in our calculations by taking the modular representation 
\[fG_0 \equiv P \mod F_0\]
Instead of solving directly for $f$, we will calculate an "auxiliary polynomial" $H$ with degree  $deg(H) < deg(F_0)$ such that
\[HG_0 \equiv 1 \mod F_0\]
in which case we can simply find $f$ by
\[f \equiv f\cdot HG_0 \equiv HP \mod F_0\]

From there, we get the updated estimates $F_1 = F_0 + f$ and $G_1 = P / F_1$.

It turns out that $H$ can be represented as another complex integral, with the following result.

\textbf{Theorem}: Given our unit splitting circle $C$ with $P = FG$ such that $F$ only has roots inside of $C$ and $G$ does not, \[HG_0 \equiv 1 \mod F_0\] is uniquely defined by \[H = 
\frac{1}{2i\pi}\oint_C \frac{1}{F_0G_0}\frac{F_0(z)-F_0(t)}{z-t}dt\]

\textit{Proof}. Let \[I(z) := \frac{1}{2i\pi}\oint_C \frac{
H(t)}{F_0(t)}\frac{F_0(z)-F_0(t)}{z-t}dt\]
We can decompose our modular equation into the form
\[\frac{1}{F_0G_0} = R + \frac{H}{F_0} + \frac{B}{G_0}\]
where $deg(A) < deg(F_0), deg(B) < deg(G_0)$ and we have that $HG_0 \equiv 1 \mod F_0$. 

We note that the map $x \rightarrow (R(x) + 
\frac{B(x)}{G_0(x)})\frac{F_0(z)-F_0(x)}{z-x}$ is analytic in the unit disk, so we can use Cauchy's integral formula to get
\[I(z) = \frac{1}{2i\pi}\oint_{|t|=1}\frac{1}{F_0(t)G_0(t)}\frac{F_0(z) -F_0(t)}{z-t}dt =\]
\[ F_0(z)\frac{1}{2i\pi}\oint_{|t|=1}\frac{H(t)}{F_0(t)}\frac{1}{z-t}dt - \frac{1}{2i\pi}\oint_{|t|=1}\frac{H(t)}{z-t}dt\]

We note that for $|z| > 1$, the second term vanishes, giving us 
\[\frac{I(z)}{F_0(z)} = \frac{1}{2i\pi}\oint_{|t|=1}\frac{H(t)}{F_0(t)}\frac{1}{z-t}dt\]

Changing the circle from $|t| = 1$ to $|t| = R$ and sending $R$ to infinity gives us a unique residue of $H(z)/F_0(z)$ and the integral is 0 for $R = \infty$, as $H$ has degree less than $F$. We conclude that $H(z) = I(z)$ for all $z$ with modulus greater than 1, so by the identity theorem, our claim is valid.

Assuming that the polynomial $F_0(z) = \phi_0 z^k + \phi_1 z^{k-1} + ... + \phi_k$, we see by linearity of integrals that $H$ would be equal to
\[H(z) = \sum_{\ell = 0}^{k-1} \left(\sum_{m=\ell+1}^k \phi_{k-m}v_{m-\ell}\right)z^\ell\]
with
\[v_m = \frac{1}{2i\pi}\oint_C \frac{t^{m-1}}{F_0G_0(t)}dt\]

Again, we can approximate this integral around the unit circle with discrete sums of values at roots of unity. We define 
\[U_m = \frac{1}{N} \sum_{j=0}^{N-1}\frac{\omega^{mj}}{P(\omega^j)}\]
One can show that $U_m$ is sufficiently close to $v_m$ given a large enough $N$ and that therefore we can make an initial guess $H_0$ of $H$ by
\[H_0(z) = \sum_{\ell = 0}^{k-1} \left(\sum_{m=\ell+1}^k \phi_{k-m}U_{m-\ell}\right)z^\ell\]
with $|H-H_0| = O(N^{n-1}e^{-\delta N})$. We then apply Newton's method to get an improved guess $H_1$ of $H$ via the following iteration:

Given $H_m$, we find $D_m$ by taking the modulus
\[H_mG_0 \equiv 1 -D_{m} \mod F_0\]
We then find $H_{m+1}$ by taking the modulus
\[H_{m+1}G_0 \equiv H_m(1 +D_{m}) \mod F_0\]
We note that if $H_0$ is a sufficiently good guess, then the sequence $\{D_m\}$ will converge quadratically to 0, as 
\[D_{m+1} \equiv 1 - H_{m+1}G_0  \equiv 1-(1-D_m)(1+D_m) \equiv D_m^2 \mod F_0\]
We iterate until reaching a desired error tolerance of $H$.

\subsection*{Discrete Fourier Transforms}
For ease of implementation, it turns out that the discrete integrals we use over the unit circle exactly correspond to Discrete Fourier transforms, which are defined as follows. Let $P$ be an $L$-dimensional vector of coefficients for our polynomial $p$, then
\[DFT(P)_m = \sum_{k=0}^{L-1} P_k \omega^{mk}\] where $\omega$ is the $L$th root of unity. Then taking a discretization of the integral over the unit circle is just $DFT(P)/L$. In addition, we can use the Fast Fourier Transform when $L$ is a power of 2 for an $O(n\log n)$ calculation of polynomial multiplication. We note that for $N = KL$ as described above:
\[W_m = \sum_{u =0}^{K-1} Y_{m,u} \omega^{(m+1)u}\]
\[U_m = \sum_{u =0}^{K-1} X_{m,u} \omega^{mu}\]
with
\[Y_{m,u} = \sum_{v=0}^{L-1}\frac{P'(\omega^{u+vK})}{P(\omega^{u+vK})}(\omega^K)^{(m+1)v}\]
\[X_{m,u} = \sum_{v=0}^{L-1}\frac{1}{P(\omega^{u+vK})}(\omega^K)^{mv}\]
where we can represent the double sums $W_m, U_m$ as sums of compositions of discrete Fourier transforms. Note that setting $K = 1$ just gives us a direct composition of DFTs. For tuning purposes, it can be helpful to have $K > 1$.

We are now ready to reveal the algorithms for estimating the factors $F(x),G(x)$ of $p(x)$ given a unit splitting circle and the number of roots inside the circle.

\subsection*{Algorithms}
We first want to take discrete Fourier transforms to estimate our integrals $W_m, U_m$:
\begin{algorithm}[H]
	\caption{DFT: Find the Fourier transforms $W_m, U_m$} 
	 \hspace*{\algorithmicindent} \textbf{Input}: $(N,k,P)$ where $P$ is a polynomial having $k < n$ roots in the unit disk, $N$ is an integer of the form $N = KL$ with $L$ is a power of $2$ such that $n < L \leq 2n$. \\
 \hspace*{\algorithmicindent} \textbf{Output}: $(W, U)$ where $W_m, U_m$ are defined for $1 \leq m \leq k$
	\begin{algorithmic}[1]
	\State Calculate $K,L$ where $N = KL$ and $L$ is a power of $2$ such that $n < L \leq 2n$.
	\State Set $W_m = 0, U_m = 0$ for $1 \leq m \leq k$.
	\State Set $\omega = e^{2i\pi/N}$ and $\omega_0 = \omega^K$.
		\For {$u=0,1,\ldots, K-1$} 
		  \State  Calculate $\alpha = FFT(P,L)$, where $\alpha_v = \sum_{j=0}^n (p_j \omega^{uj})\omega_0^{vj}$ for $v \in 0:L-1$.
		  \State  Calculate $\beta = FFT(P',L)$, where $\beta_v = \sum_{j=0}^{n-1} ((j+1)p_{j+1} \omega^{uj})\omega_0^{vj}$ for $v \in 0:L-1$.
		  \State Calculate $\gamma_v = 1/\alpha_v$ for $v \in 0:L-1$.
		  \State Calculate $x = FFT(\gamma, L)$, where $x_m = \sum_{v=0}^{L-1} \gamma_v \omega_0^{mv}$ for $m \in 0:L-1$.
		  \State Set $U_m = U_m + x_m \omega^{mu}$ for $m \in 0:k-1$.
		  \State Calculate $\delta_v = \beta_v/\alpha_v$ for $v \in 0:L-1$.
		  \State Calculate $y = FFT(\delta, L)$, where $y_m = \sum_{v=0}^{L-1} \delta_v \omega_0^{m v}$ for $m \in 0:L-1$.
		  \State Set $W_m = W_m + y_{m+1} \omega^{(m+1) u}$ for $m \in 0:k-1$.
		    \EndFor
    \State Set $W_m = W_m/N$ and $U_m = U_m/N$.
	\end{algorithmic} 
\end{algorithm}

Now that we have the appropriate integrals, we would like to plug them in to get initial guesses for our factor $F$ and auxiliary polynomial $H$. 
\begin{algorithm}[H]
	\caption{RES: Initial approximation of $F,H$ given the respective contour integrals} 
	 \hspace*{\algorithmicindent} \textbf{Input}: $(N,k,P,\delta, \varepsilon)$ where $P$ is a polynomial having $k < n$ roots in the unit disk with no roots in the annulus $(e^{-\delta}, e^{\delta})$, $N$ is an integer of the form $N = KL$ with $L$ is a power of $2$ such that $n < L \leq 2n$, $\varepsilon$ is a preset error tolerance. \\
 \hspace*{\algorithmicindent} \textbf{Output}: $(F_0, H_0)$ which are approximations of $F,H$ respectively.
	\begin{algorithmic}[1]
	\State Calculate the contour integrals $(W, U)$ approximately using DFT($P,k,N)$.
	\State Approximate the factors $\phi_m$ in $F_0 = z^k + \phi_1z^{k-1} + ... \phi_k$ using the recursive formulas
	\[\phi_1 = -W_i\]
	\[\phi_m = -\frac{1}{m}(W_1\phi_{m-1} + ... + W_{m-1}\phi_1 + W_m)\]
    \State Set the polynomial $H_0(z) = \sum_{\ell = 0}^{k-1} \left(\sum_{m=\ell+1}^k \phi_{k-m}U_{m-\ell}\right)z^\ell$ and keep the coefficients as a $k$-dimensional vector.
	\end{algorithmic} 
\end{algorithm}

We now improve $H_0$ for use in Newton's method:

\begin{algorithm}[H]
	\caption{AUX: Improves approximation of auxiliary polynomial $H$} 
	 \hspace*{\algorithmicindent} \textbf{Input}: $(F_0, G_0, H_0, \varepsilon)$ where $F_0,G_0,H_0$ are approximations of $F,G,H$, $\varepsilon$ is a preset error tolerance. \\
 \hspace*{\algorithmicindent} \textbf{Output}: $(H_1)$, which is an approximation of $H$ such that $H_1G_0 \equiv 1 -D \mod F_0$ with $|D| < \varepsilon$
	\begin{algorithmic}[1]
	\State Set $H_1 := H_0$.
	\For{$k=0,1,2,\ldots$}
	\State Find $D$ by the modular relation $H_1G_0 \equiv 1 -D \mod F_0$. 
	\If{$||D||_1 < \varepsilon$}
	\State Return $H_1$.
	\ElsIf{$||D||_1 > 1$}
	\State Return 0.
	\EndIf
	\State Set $H_{1} \equiv H_1(1+D) \mod F_0$
	\EndFor
	\end{algorithmic} 
\end{algorithm}

We now plug in our improved estimate of $H$ to the Newton-Schonhage method:
\begin{algorithm}[H]
	\caption{NS: Improves approximation of factor $F$} 
	 \hspace*{\algorithmicindent} \textbf{Input}: $(P, F_0, H_0, \varepsilon)$ where $P$ is a polynomial, $F_0, H_0$ are approximations of $F,H$, $\varepsilon$ is a preset error tolerance. \\
 \hspace*{\algorithmicindent} \textbf{Output}: $(F_1,G_1)$, which are approximations of $F,G$ such that $|P-F_1G_1| < \varepsilon|P|$.
	\begin{algorithmic}[1]
	\State Set $F_1 := F_0$
	\State Set $G_1 := P/F_1$ via Euclidean division.
	\State Calculate $\varepsilon_0 = |P-F_1G_1|/|P|$.
	\For{k = 0,1,2....}
	\If{$\varepsilon_0 < \varepsilon$}
	\State Return $(F_1,G_1)$.
	\ElsIf{$\varepsilon_0 > 1$}
	\State Return (0,0).
	\EndIf
	\State Set $H := $ AUX$(F_1,G_1,H_1,\varepsilon_0)$.
	\State Calculate $f := H_1P \mod F_1$.
	\State Set $F_1 := F_1 + f$, $G_1 = P/F_1$.
	\EndFor
	\end{algorithmic} 
\end{algorithm}
Note that the above algorithm could fail to converge if the approximation $F_0$ of $F$ is too far off. We combine all of our steps for the final algorithm:
\begin{algorithm}[H]
	\caption{FCS: Factorization given a splitting circle} 
	 \hspace*{\algorithmicindent} \textbf{Input}: $(P, k, \delta, \varepsilon)$ where $P$ is a polynomial of degree $n \geq 2$ possessing $k< n$ roots in the unit disk and none in the annulus $(e^{-\delta}, e^{\delta})$, $\varepsilon > 0$ being the preset error tolerance \\
 \hspace*{\algorithmicindent} \textbf{Output}: $(F_1,G_1)$, which are approximations of $F,G$ such that $|P-F_1G_1| < \varepsilon|P|$.
	\begin{algorithmic}[1]
	\State Set $L$ to be the power of 2 such that $n < L \leq 2n$
	\State Set $K = \max(1/2\delta, 2)$ and set $N = KL$.
	\For{k = 0,1,2....}
	\State Find the initial guesses $(F_0,G_0) = RES(P,k,\delta, N)$
	\State Find the updated factors $F_1,G_1 = NS(P,F_0,H_0,\varepsilon)$
    \If{$(F_1,G_1) = (0,0)$}
	\State Set $N := 2N$.
	\Else 
	\State Return $(F_1,G_1)$. 
	\EndIf
	\EndFor
	\end{algorithmic} 
\end{algorithm}

\section*{Complete Factorization}
Now that we are capable of finding a splitting circle and using it to factor our degree $n$ polynomial $P$ into two nontrivial factors $F,G$, we can use our routine $CTR0$ repeatedly until we've reduced the problem to linear factors $L_1,L_2,...L_n$, with
\[|P - L_1L_2...L_n| < \varepsilon |P|\]
This would occur if each intermediary call to $CTR0$ left us with $k$ factors $P_1,...,P_k$ such that 
\[|P - P_1P_2...P_k| < \frac{k}{n}\varepsilon |P|\]
While $k < n$, we can factorize one of the polynomials $P_i$ with degree $>1$. WLOG, let this be $P_1$, with a factorization of 
\[|P-FG| < \varepsilon_k|P_1|\]
which would give
\[|P - FGP_2...P_k| < \frac{k}{n}\varepsilon |P| + \varepsilon_k|P_1|\cdot|P_2...P_k|\]
This would be upper bounded by $\frac{k+1}{n}\cdot \varepsilon|P|$ if $\varepsilon_k < \frac{\varepsilon}{n}\frac{|P|}{n|P_1|\cdot|P_2...P_k|}$.
From Schonhage, we have the inequality
\[|Q|\cdot|R| < 2^{\deg(QR)-1}|QR|\]
which implies
\[|P_1|\cdot|P_2...P_k| < 2^{n-1}|P_1...P_k| < 2^{n-1}(1+\frac{k}{n}\varepsilon)|P| < 2^n|P|\]
so our constraint is modified if we choose 
\[\varepsilon_k = 2^{-n}\frac{\varepsilon}{n}\]
In summary, if each factorization step of the algorithm $CTR0$ uses a precision of $2^{-n}\varepsilon/n$, then we have our approximate linear factorization. That gives us the following algorithm:
\begin{algorithm}[H]
	\caption{FACT: Approximate complete factorization} 
	 \hspace*{\algorithmicindent} \textbf{Input}: $(P, \varepsilon)$ where $P$ is a polynomial of degree $n \geq 1$, $0< \varepsilon <1 $ being the preset error tolerance \\
 \hspace*{\algorithmicindent} \textbf{Output}: $(L_1,...,L_n)$, which are linear factors such that such that $|P-L_1L_2...L_n| < \varepsilon|P|$.
	\begin{algorithmic}[1]
	\State If $n=1$, return $P$.
	\State Otherwise, we calculate $(F,G) = CTR0(P,\frac{2^{-n}\varepsilon}{n})$.
	\State Let $k = deg(F)$. Then we calculate $(L_1,...,L_k) = FACT(F,\varepsilon)$ and $(L_{k+1},...,L_n) = FACT(G,\varepsilon)$
	\State Return $(L_1,...,L_n)$.
	\end{algorithmic} 
\end{algorithm}

\section*{References}

1. Gourdon, Xavier. \textit{Combinatoire, Algorithmique et Géométrique des Polynômes.} PhD Thesis, École Polytechnique, Paris, 1996.

2. Schonhage, A. Equation solving in terms of computational complexity. In \textit{Proceedings of the International Congress of Mathematicians} (1987).

3. Schatzle, R. Diploma Thesis
 Universitat Bonn, 1996.

\end{document}